# An Extension of Tychonoff's Fixed Point Theorem in Hausdorff Quasi-locally Convex Topological Vector Spaces

Jinlu Li

Department of Mathematics

Shawnee State University

Portsmouth, Ohio 45662

USA

jli@shawnee.edu

**Abstract**

In this paper, we introduce the concepts of weaknorm, quasi-weaknorm on real vector spaces. By these concepts, we introduce the concept of quasi-locally convex topological vector spaces, which include locally convex topological vector spaces as special cases. By the Fan-KKM theorem, we prove a fixed point theorem in quasi-locally convex topological vector spaces, that is a natural extension of Tychonoff's fixed point theorem in locally convex topological vector spaces. Then we provide an example to show that this extension is a proper extension.

1. Introduction

In 1912, Brouwer proved the first fixed-point theorem, which is in Euclidean spaces.

**Brouwer's fixed-point theorem** [1]. Every continuous mapping from the $n$-simplex to itself has a fixed point.

In 1930, Schauder extended Brouwer's fixed-point theorem from Euclidean spaces to Banach spaces.

**Schauder's fixed-point theorem** [9]**.** Let $C$ be a nonempty closed convex subset of a Banach space $X$. If $f : C \to C$ is continuous with a compact image, then $f$ has a fixed point.

Meanwhile, Schauder had the well-known conjecture, which has been named as Schauder's Conjecture: *Every continuous function, from a nonempty compact and convex set in a (Hausdorff) topological vector space into itself, has a fixed point.* [Problem 54 in The Scottish Book].

R. Cauty in [4] proposed an answer to the Schauder's Conjecture. In the international conference of Fixed Point Theory and its Applications in 2005, T. Dobrowolski remarked that there is a gap in the proof. Therefore, Schauder's Conjecture is still unsolved.

To generalize the underlying spaces in fixed point theory, in 1934, Tychonoff extended Schauder's fixed point theorem from Banach spaces to locally convex topological vector space.

**Tychonoff's fixed point theorem** [12]**:** Let $X$ be a Hausdorff locally convex topological vector space. For any nonempty compact convex set $C$ in $X$, any continuous function $f : C \to C$ has a fixed point.

Schauder's fixed-point theorem and Tychonoff's fixed point theorem have been extensively applied in many fields of mathematics. Many authors have extended these theorems regarding to the considered mappings and the underlying spaces. These theorems have been proved by using different ways (see [2], [5−7], [11]).

The Fan-KKM theorem has played an important role in nonlinear analysis, fixed point theory, optimization theory, variational analysis, etc. In [5], Olga Hadžić and Endre Pap used the Fan-KKM theorem proved the Tychonoff's fixed point theorem, the Brouwer's fixed-point theorem and the Schauder's fixed-point theorem. In [8], Park proved some fixed point theorems by the Fan-KKM theorem. It shows that the Fan-KKM theorem is indeed a very powerful tool to prove some fixed point theorems and to prove the solvability of some global optimization problems.

This paper is organized as follows: In section 2, we briefly recall the concepts of locally convex topological vector spaces; in section 3, we introduce the concepts of weaknorm and quasi-weaknorm on vector spaces that induce the concept of quasi-locally convex topological vector spaces. Some properties are provided including that the concept of quasi-locally convex is an extension of locally convex; in section 4, we prove a fixed point theorem in Hausdorff quasi-locally convex topological vector spaces by using the Fan-KKM theorem; in section 5, we give an example to show that the concept of quasi-locally convex is a proper extension of locally convex. It implies that our fixed point theorem proved in section 4 is a proper extension of the Tychonoff's fixed point theorem.

2.  **Preliminaries**

In this section, we briefly recall the definitions of locally convex topological vector spaces.

2.1. Locally convex topological vector spaces defined via convex sets

Let $X$ be a real vector space. A subset $C$ in $X$ is called

1. Convex if for all $x$, $y$ in $C$, and $0 \leq t \leq 1$, $tx + (1 - t)y$ is in $C$.
2. Balanced if for all $x$ in $C$ and for $|\alpha| \leq 1$, $\alpha x$ is in $C$.

3. Absorbing if the union of $tC$ over all $t \in \mathbb{R}$ is all of $X$.

**Definition (first version).** A topological vector space is said to be locally convex if the origin has a local base of balanced, convex and absorbing sets.

Because translation is (by definition of "topological vector space") continuous, all translations are homeomorphisms, so every base for the neighborhoods of the origin can be translated to a base for the neighborhoods of any given vector.

**Definition (second version).** If a topological vector space is locally convex, then there is a base at the origin consisting of convex open sets.

2.2 Locally convex topological vector spaces defined via seminorms

A mapping $p: X \to \mathbb{R}^+$ is called a seminorm on $X$ if it satisfies the following conditions:

$S_1$. $p(x) \geq 0$, for all $x \in X$;
$S_2$. $p(\alpha x) = |\alpha| p(x)$, for every scaler $\alpha$ and for all $x \in X$;
$S_3$. $p(x + y) \leq p(x) + p(y)$, for all $x, y \in X$.

**Definition (third version)** A topological vector space $(X, \tau)$ is a locally convex topological vector space, if and only if there is a family of seminorms $\{p_\lambda\}_{\lambda \in \Lambda}$ on $X$ such that the initial topology $\tau$ on $X$ is induced by the seminorms.

From the definition, $\tau$ is the coarsest topology on $X$, for which, the following mappings

$$p_{\lambda,y}(\cdot) = p_\lambda(\cdot - y): X \to \mathbb{R}^+, \text{ for every } y \in X \text{ and for every } \lambda \in \Lambda, \tag{2.1}$$

are $\tau$-continuous. For every $y \in X$, a base of neighborhoods of $y$ for this topology is obtained as follows:

$$B_{\Gamma,\varepsilon}(y) = \{x \in X: p_\lambda(x - y) < \varepsilon, \text{ for every } \lambda \in \Gamma\}, \tag{2.2}$$

for every finite subset $\Gamma$ of $\Lambda$, and for every $\varepsilon > 0$.

2.3. The Fan-KKM theorem

For easy reference, we recall the KKM mappings and the Fan-KKM theorem below. KKM mappings are defined in vector spaces and the underlying spaces of the Fan-KKM Theorem are Hausdoff topological vector spaces. For more details, the readers are referred to Fan [5] and Park [8−9].

Let $C$ be a nonempty subset of a vector space $X$. A set-valued mapping $F: C \to 2^X \setminus \{\emptyset\}$ is called a KKM mapping if, for any finite subset $\{x_1, x_2, \ldots, x_k\}$ of $C$, we have

$$\text{co}\{x_1, x_2, \ldots, x_k\} \subseteq \bigcup_{1 \leq i \leq k} F(x_i),$$

where $\text{co}\{x_1, x_2, \ldots, x_k\}$ denotes the convex hull of $\{x_1, x_2, \ldots, x_k\}$.

**Fan-KKM Theorem**. *Let C be a nonempty closed convex subset of a Housdorff topological vector space and let $F: C \to 2^C \setminus \{\emptyset\}$ be a KKM mapping with closed values. If there exists a point $x_0 \in C$ such that $F(x_0)$ is a compact subset of C, then*

$$\bigcap_{x \in C} F(x) \neq \emptyset.$$

### 3. Quasi-locally convex Topological Vector Spaces

In this paper, we only consider real vector spaces. In this section, we introduce new concepts of weaknorm, quasi-weaknorm on real vector spaces, which induce the concept of quasi-locally convex topological vector spaces. We also give some properties of weaknorms.

**Definition 3.1.** Let $X$ be a vector space. A mapping $p: X \to \mathbb{R}^+$ is called a weaknorm on $X$ if it satisfies the following conditions:

$W_1$. $p(x) \geq 0$, for all $x \in X$ and $p(\theta) = 0$;
$W_2$. $p(-x) = p(x)$, for all $x \in X$;
$W_3$. For any elements $x_1, x_2$ of $X$, and $0 \leq \alpha \leq 1$, one has

$$p(\alpha x_1 + (1 - \alpha) x_2) \leq \alpha p(x_1) + (1 - \alpha) p(x_2).$$

It is clear to see that the condition $W_3$ implies that for any finite set of distinct elements $x_1, x_2, \ldots, x_n$ of $X$, and positive numbers $\alpha_1, \alpha_2, \ldots, \alpha_n$ satisfying $\sum_{i=1}^n \alpha_i = 1$, one has

$$p(\sum_{i=1}^n \alpha_i x_i) \leq \sum_{i=1}^n \alpha_i p(x_i). \tag{3.1}$$

The vector space $X$ equipped with a weaknorm $p$ is called a weaknormed vector space, denoted by $(X, p)$.

**Lemma 3.2.** Let $X$ be a vector space. Every seminorm on $X$ is a weaknorm on this vector space. Therefore, every seminormed vector space is a weaknormed vector space.

**Definition 3.3.** Let $X$ be a vector space. A mapping $q: X \to \mathbb{R}^+$ is called a quasi-weaknorm on $X$ if there are a weaknorm $p$ on $X$ and a strictly increasing continuous function $\varphi: \mathbb{R}^+ \to \mathbb{R}^+$ such that

$W_4$. $q(x) \leq \varphi(p(x))$, for all $x \in X$;
$W_5$. $\varphi(0) = 0$.

The vector space $X$ equipped with a quasi-weaknorm $q$ is called a quasi-weaknormed vector space, denoted by $(X, q)$, in which $p$ is called the adjoint weaknorm and $\varphi$ is called the associated weighted function of $q$.

**Lemma 3.4.** *Let X be a vector space. Every weaknorm on X is a quasi-weaknorm on this vector space and every weaknormed vector space is a quasi-weaknormed vector space.*

*Proof.* Take $\varphi(t) = t$, for every $t \in \mathbb{R}^+$. □

Therefore, for any vector space $X$, every seminorm on $X$ is a quasi-weaknorm on this vector space and every seminormed vector space is a quasi-weaknormed vector space.

**Definition 3.5**. Let $(X, \tau)$ be a topological vector space. If there is a family of $\tau$-continuous quasi-weaknorms $\{q_\lambda\}_{\lambda \in \Lambda}$ on $X$ with $\tau$-continuous adjoint weaknorms $\{p_\lambda\}_{\lambda \in \Lambda}$ and associated weighted functions $\{\varphi_\lambda\}_{\lambda \in \Lambda}$ such that the initial topology $\tau$ on $X$ is induced by the quasi-weaknorms $\{q_\lambda\}_{\lambda \in \Lambda}$, then $(X, \tau)$ is called a quasi-locally convex topological vector space,

From the definition, $\tau$ is the coarsest topology on $X$, for which, the following mappings

$$q_\lambda(\cdot - y): X \to \mathbb{R}^+, \text{ for any } y \in X \text{ and for every } \lambda \in \Lambda, \tag{3.2}$$

are $\tau$-continuous. Furthermore, for any $y \in X$ and for every $\lambda \in \Lambda$, the following mappings

$$p_\lambda(\cdot - y): X \to \mathbb{R}^+ \quad \text{and} \quad \varphi_\lambda(p_\lambda(\cdot - y)): X \to \mathbb{R}^+ \tag{3.3}$$

are both $\tau$-continuous.

**Definition 3.6**. A family of quasi-weaknorms $\{q_\lambda\}_{\lambda \in \Lambda}$ on a topological vector space $(X, \tau)$ is called total whenever for some $x \in X$, $q_\lambda(x) = 0$ holds for every $\lambda \in \Lambda$, then it is necessary to have $x = \theta$.

**Lemma 3.7.** Every locally convex topological vector space is a quasi-locally convex topological vector space.

*Proof.* Let $(X, \tau)$ be a locally convex topological vector space, in which the initial topology $\tau$ on $X$ is induced by a family of seminorms $\{p_\lambda\}_{\lambda \in \Lambda}$. By definition of locally convex topological vector spaces and (2.1), one has that, for every $y \in X$ and for every $\lambda \in \Lambda$, the mapping

$$p_\lambda(\cdot - y): X \to \mathbb{R}^+$$

is $\tau$-continuous. Then for every $\lambda \in \Lambda$, we define

$$\varphi_\lambda(t) = t, \text{ for every } t \in \mathbb{R}^+$$

and

$$q_\lambda(x) = \varphi_\lambda(p_\lambda(x)) = p_\lambda(x), \text{ for } x \in X.$$

This lemma follows immediately. □

**Remarks 3.8**. The reasons that a mapping $p: X \to \mathbb{R}^+$ in Definition 3.1 is named as a weaknorm are that $p$ is not required to satisfy the following absolute homogeneity and triangle inequality:

$S_2$. $p(\lambda x) = |\lambda| p(x)$, for every scalar $\lambda$ and for all $x \in X$,
$S_3$. $p(x + y) \leq p(x) + p(y)$, for all $x, y \in X$.

We give some properties of weaknorms in a lemma below.

**Lemma 3.9**. Let $(X, \tau)$ be a topological vector space equipped with a weaknorm $p: X \to \mathbb{R}^+$. For any $x \in X$, one has that

(i). If $0 < \alpha \leq 1$, then

$$p(\alpha x) \leq \alpha p(x),$$

(ii). If $\alpha > 1$, then

$$p(\alpha x) \geq \alpha p(x).$$

*Proof.* By (3.1), for $0 < \alpha \leq 1$, we have

$$p(\alpha x) = p(\alpha x + (1-\alpha)\theta) \leq \alpha p(x) + (1-\alpha)p(\theta) = \alpha p(x).$$

If $\alpha > 1$, we have

$$p(x) = p(\tfrac{1}{\alpha}\alpha x) = p(\tfrac{1}{\alpha}\alpha x + (1-\tfrac{1}{\alpha})\theta) \leq \tfrac{1}{\alpha}p(\alpha x) + (1-\tfrac{1}{\alpha})p(\theta) = \tfrac{1}{\alpha}p(\alpha x).$$

It implies

$$p(\alpha x) \geq \alpha p(x).$$

## 4. A fixed point theorem in Housforrf quasi-locally convex topological vector spaces

**Theorem 4.1.** *Let $(X, \tau)$ be a Housforrf quasi-locally convex topological vector space with a total family of quasi-weaknorms. Let C be a nonempty compact convex subset of X. Then every continuous mapping from C to itself has a fixed point.*

*Proof.* Let $\{q_\lambda\}_{\lambda \in \Lambda}$ be a total family of $\tau$-continuous quasi-weaknorms on the topological vector space $X$ which defines (induces) the topology $\tau$, in which $\{q_\lambda\}_{\lambda \in \Lambda}$ has a family of $\tau$-continuous adjoint weaknorms $\{p_\lambda\}_{\lambda \in \Lambda}$ and associated with the family of weighted functions $\{\varphi_\lambda\}_{\lambda \in \Lambda}$.

Let $f: C \to C$ be a continuous mapping. The totality of $\{q_\lambda\}_{\lambda \in \Lambda}$ implies that a point $x_0 \in C$ is a fixed point of the mapping $f$ if and only if

$$q_\lambda(x_0 - f(x_0)) = 0, \text{ for every } \lambda \in \Lambda.$$

It follows that the mapping $f$ has a fixed point if and only if

$$\bigcap_{\lambda \in \Lambda}\{x \in C : q_\lambda(x - f(x)) = 0\} \neq \emptyset. \tag{4.1}$$

From the compactness of $C$, to prove (4.1), we only need to prove that the following family has the finite intersection property

$$\Big\{\{x \in C : q_\lambda(x - f(x)) = 0\} : \lambda \in \Lambda\Big\}.$$

Let $m$ be an arbitrary positive integer and let $\{\lambda_1, \lambda_2, \ldots, \lambda_m\}$ be an arbitrary finite subset of $\Lambda$. To prove the finite intersection property, therefore, we prove that

$$\bigcap_{1 \leq k \leq m}\{x \in C : q_{\lambda_k}(x - f(x)) = 0\} \neq \emptyset. \tag{4.2}$$

We first prove that, for any $\delta > 0$,

$$\bigcap_{1 \leq k \leq m}\{x \in C : q_{\lambda_k}(x - f(x)) < \delta\} \neq \emptyset. \tag{4.3}$$

From condition $W_5$, we have

$$q_{\lambda_k}(x - f(x)) \leq \varphi_{\lambda_k}\left(p_{\lambda_k}(x - f(x))\right), \text{ for } k = 1, 2, \ldots, m \text{ and for every } x \in C.$$

It implies that

$$\bigcap_{1 \leq k \leq m}\left\{x \in C: \varphi_{\lambda_k}\left(p_{\lambda_k}(x - f(x))\right) < \delta\right\} \subseteq \bigcap_{1 \leq k \leq m}\{x \in C: q_{\lambda_k}(x - f(x)) < \delta\}. \quad (4.4)$$

Therefore, to prove (4.3), it is enough to show that

$$\bigcap_{1 \leq k \leq m}\left\{x \in C: \varphi_{\lambda_k}\left(p_{\lambda_k}(x - f(x))\right) < \delta\right\} \neq \emptyset. \quad (4.5)$$

Notice that, for $k = 1, 2, \ldots, m$, every $\varphi_{\lambda_k}: \mathbb{R}^+ \to \mathbb{R}^+$ is a strictly increasing continuous function and $\varphi_{\lambda_k}(0) = 0$. Then $\varphi_{\lambda_k}^{-1}(\delta) > 0$ and (4.5) is equivalent to

$$\bigcap_{1 \leq k \leq m}\{x \in C: p_{\lambda_k}(x - f(x)) < \varphi_{\lambda_k}^{-1}(\delta)\} \neq \emptyset. \quad (4.6)$$

Let $\beta = \min\{\varphi_{\lambda_k}^{-1}(\delta): k = 1, 2, \ldots, m\}$. Then $\beta > 0$ and it follows that

$$\{x \in C: \sum_{k=1}^{m} p_{\lambda_k}(x - f(x)) < \beta\}$$

$$\subseteq \bigcap_{1 \leq k \leq m}\{x \in C: p_{\lambda_k}(x - f(x)) < \beta\}.$$

$$\subseteq \bigcap_{1 \leq k \leq m}\{x \in C: p_{\lambda_k}(x - f(x)) < \varphi_{\lambda_k}^{-1}(\delta)\}. \quad (4.7)$$

Next, we use the Fan-KKM Theorem to show that

$$\{x \in C: \sum_{k=1}^{m} p_{\lambda_k}(x - f(x)) < \beta\} \neq \emptyset. \quad (4.8)$$

For $\beta > 0$ as given above, assume, on the contrary, that (4.8) does not hold; that is,

$$\sum_{k=1}^{m} p_{\lambda_k}(x - f(x)) \geq \beta, \text{ for every } x \in C. \quad (4.9)$$

Based on the mapping $f$, we define a set-valued mapping $F: C \to 2^C \setminus \{\emptyset\}$ as follows:

$$F(x) = \{z \in C: \sum_{k=1}^{m} p_{\lambda_k}(x - f(z)) \geq \beta\}, \text{ for } x \in C.$$

From the hypothesis (4.9), we see that $x \in F(x)$, and therefore, $F(x) \neq \emptyset$, for every $x \in C$. The $\tau$-continuity of $f: C \to C$ and the $\tau$-continuouity adjoint weaknorms $\{p_\lambda\}_{\lambda \in \Lambda}$ (see Definition 3.5 and (3.3)) imply that $F(x)$ is $\tau$-closed, for every $x \in C$. Next, we show that the set-valued mapping $F: C \to 2^C \setminus \{\emptyset\}$ is a KKM mapping.

For any given positive integer $n$, take arbitrary $n$ distinct points $x_1, x_2, \ldots, x_n \in C$. For any positive numbers $t_1, t_2, \ldots, t_n$ satisfying $\sum_{i=1}^{n} t_i = 1$, let $y = \sum_{i=1}^{n} t_i x_i$. We show that

$$y \in \bigcup_{1 \leq i \leq n} F(x_i). \quad (4.10)$$

Assume, by the way of contradiction, that (4.10) does not hold. Then we have

$$y \notin F(x_i), \text{ for every } i = 1, 2, \ldots, n.$$

It is

$$\sum_{k=1}^m p_{\lambda_k}(x_i - f(y)) < \beta, \text{ for all } i = 1, 2, \ldots, n. \tag{4.11}$$

From the assumption (4.9), the inequality (4.11), and by the condition W$_3$ of the weaknorm $p_{\lambda_k}$, it follows that

$$\begin{aligned}
\beta &\leq \sum_{k=1}^m p_{\lambda_k}(y - f(y)) \\
&= \sum_{k=1}^m p_{\lambda_k}(\sum_{i=1}^n t_i x_i - f(y)) \\
&= \sum_{k=1}^m p_{\lambda_k}(\sum_{i=1}^n t_i(x_i - f(y))) \\
&\leq \sum_{k=1}^m \sum_{i=1}^n t_i p_{\lambda_k}(x_i - f(y)) \\
&= \sum_{i=1}^n t_i \sum_{k=1}^m p_{\lambda_k}(x_i - f(y)) \\
&< \sum_{i=1}^n t_i \beta \\
&= \beta.
\end{aligned}$$

It is a contradiction. It implies that $F: C \to 2^C \setminus \{\emptyset\}$ is a KKM mapping with nonempty closed values. Since $C$ is compact, from the Fan-KKM Theorem, we have

$$\bigcap_{x \in C} \{z \in C : \sum_{k=1}^m p_{\lambda_k}(x - f(z)) \geq \beta\} = \bigcap_{x \in C} F(x) \neq \emptyset.$$

Then, there is $z_0 \in C$ satisfying

$$\sum_{k=1}^m p_{\lambda_k}(x - f(z_0)) \geq \beta, \text{ for every } x \in C.$$

In particularly, if we take $x = f(z_0)) \in C$, we get

$$0 = \sum_{k=1}^m p_{\lambda_k}(f(z_0) - f(z_0)) \geq \beta.$$

It is a contradiction. So (4.8) is proved. Then, (4.6) follows from (4.7) and (4.8). Meanwhile, (4.6) proves (4.5) because they are equivalent. Then (4.3) is proved by (4.5) and (4.4).

Hence, we proved that, for any $\delta > 0$,

$$\bigcap_{1 \leq k \leq m} \{x \in C : q_{\lambda_k}(x - f(x)) < \delta\} \neq \emptyset. \tag{4.3}$$

From the continuity of $f$ and the continuity of every $q_{\lambda_k}$, (4.3) implies that

$$\bigcap_{1 \leq k \leq m} \{x \in C : q_{\lambda_k}(x - f(x)) \leq \delta\}$$

is a nonempty $\tau$-closed subset of $C$. Take a strictly decreasing sequence of positive numbers $\{\delta_u\}$ with limit 0; i.e., $\delta_u \downarrow 0$, as $u \to \infty$. Then, for $w > v$, we have

$$\bigcap_{1 \leq k \leq m} \{x \in C : q_{\lambda_k}(x - f(x)) \leq \delta_w\} \subseteq \bigcap_{1 \leq k \leq m} \{x \in C : q_{\lambda_k}(x - f(x)) \leq \delta_v\}.$$

That is, $\{\bigcap_{1 \leq k \leq m} \{x \in C : q_{\lambda_k}(x - f(x)) \leq \delta_u\}\}$ is a decreasing (with respect to inclusions) sequence of nonempty closed subsets of the compact set $C$. It follows that

$$\bigcap_{1 \leq k \leq m} \{x \in C : q_{\lambda_k}(x - f(x)) = 0\} = \bigcap_{u=1}^\infty \bigcap_{1 \leq k \leq m} \{x \in C : q_{\lambda_k}(x - f(x)) \leq \delta_u\} \neq \emptyset.$$

Hence (4.2) is proved. Since $C$ is compact, by the finite intersection property, (4.1) follows immediately. □

## 5. An example: quasi-locally convexity is a proper extension of local convexity

From Lemma 3.7, we see that Fixed Point Theorem 4.1 in Hausdorff quasi-locally convex topological vector spaces is an extension of Tychonoff fixed point theorem in Hausdorff locally convex topological vector spaces. One more step further, in this section, we provide a concrete example of topological vector space, which is quasi-locally convex but not locally convex. Therefore, it shows that the quasi-locally convexity of topological vector spaces is a proper extension of the concept of local convexity of topological vector spaces. Consequently, it implies that Fixed Point Theorem 4.1 proved in the previous section properly extends Tychonoff fixed point theorem.

**Example 5. 1**. Let $S$ denote the set of real sequences. For every given $r \in (0, 1)$, we define a subspace $l_r$ of $S$ as below

$$l_r = \{\{t_i\} \in S : \sum_{i=1}^{\infty} |t_i|^r < \infty\}.$$

The space $l_r$ is a topological vector space equipped with the topology $\tau_q$ defined by the following functional:

$$q(\{t_i\}) = \sum_{i=1}^{\infty} |t_i|^r, \text{ for every } \{t_i\} \in S.$$

Then, the topological vector space $(l_r, \tau_q)$ is quasi-locally convex. However, it is known that $(l_r, \tau_q)$ is not locally convex.

*Proof.* We give a proof of the statement that $(l_r, \tau_q)$ is not locally convex in the Appendix. Here, we prove that $(l_r, \tau_q)$ is a quasi-locally convex topological vector space. For any positive integer $\lambda$, we define a mapping $q_\lambda : l_r \to R^+$ as follows

$$q_\lambda(\{t_i\}) = \sum_{i=1}^{\lambda} |t_i|^r, \text{ for every } \{t_i\} \in l_r. \tag{5.1}$$

It is clear that $q_\lambda : l_r \to \mathbb{R}^+$ is $\tau_q$-continuous. Next, we show that, for every $\lambda = 1, 2, \ldots$, the mapping $q_\lambda : l_r \to \mathbb{R}^+$ is a quasi-weaknorm on $l_r$. To this end, for every $\lambda = 1, 2, \ldots$, we define a mapping $p_\lambda : l_r \to \mathbb{R}^+$ by

$$p_\lambda(\{t_i\}) = \sum_{i=1}^{\lambda} |t_i|, \text{ for every } \{t_i\} \in l_r. \tag{5.2}$$

We next show that $p_\lambda : l_r \to \mathbb{R}^+$ is $\tau$-continuous. Suppose that, for $j = 1, 2, \ldots, \{t_{ij}\} \in l_r$ and $\{y_i\} \in l_r$ satisfy that $\{t_{ij}\} \xrightarrow{\tau_q} \{y_i\}$, as $j \to \infty$. That is,

$$q_\lambda(\{t_{ij} - y_i\}) = \sum_{i=1}^{\infty} |t_{ij} - y_i|^r \to 0, \text{ as } j \to \infty.$$

Then, for every fixed positive integer $\lambda = 1, 2, \ldots$, we have

$$\sum_{i=1}^{\lambda} |t_{ij} - y_i|^r \to 0, \text{ as } j \to \infty.$$

For every given positive integer $\lambda$, it implies that

$$p_\lambda(\{t_{ij} - y_i\}) = \sum_{i=1}^{\lambda} |t_{ij} - y_i| \to 0, \text{ as } j \to \infty.$$

Applying the following triangle inequality

$$(x + y)^r \leq x^r + y^r, \text{ for } x, y \geq 0, \text{ for } 0 < r < 1,$$

We obtain that, for every fixed positive integer $\lambda = 1, 2, \ldots,$ $p_\lambda$ is a $\tau_q$-continuous weaknorm on $l_r$ (as a matter of fact, it is a seminorm). On the other hand, from (5.1), for every $\{t_i\} \in l_r$, we have

$$q_\lambda(\{t_i\}) = \sum_{i=1}^{\lambda} |t_i|^r \leq \left(\sum_{i=1}^{\lambda} |t_i|\right)^r \lambda^{1-r}. \tag{5.3}$$

For every $\lambda = 1, 2, \ldots,$ we define a function $\varphi_\lambda: \mathbb{R}^+ \to \mathbb{R}^+$ by

$$\varphi_\lambda(u) = \lambda^{1-r} u^r, \text{ for } u \in \mathbb{R}^+. \tag{5.4}$$

Then, $\varphi_\lambda: \mathbb{R}^+ \to \mathbb{R}^+$ is continuous, strictly increasing. By (5.2), (5.3) and (5.4), it implies that

$$q_\lambda(\{t_i\}) \leq \varphi_\lambda(p_\lambda(\{t_i\})), \text{ for every } \{t_i\} \in l_r.$$

Hence it satisfies conditions $W_5$ and $W_6$. Then $q_\lambda: l_r \to \mathbb{R}^+$ is a $\tau_q$-continuous quasi-weaknorm with the $\tau_q$-continuous adjoint weaknorms $p_\lambda$ defined in (5.2) and associated with the weighted functions $\varphi_\lambda$ defined in (5.4).

Notice that

$$q_\lambda(\{t_i\}) \leq q(\{t_i\}), \text{ for } \lambda = 1, 2, \ldots, \tag{5.5}$$

and

$$q(\{t_i\}) = \lim_{\lambda \to \infty} q_\lambda(\{t_i\}), \text{ for every } \{t_i\} \in l_r.$$

One can see that the topology on $l_r$ defined by the mapping $q$ is induced by the family of the quasi-weaknorms $\{q_\lambda: \lambda = 1, 2, \ldots\}$ on $l_r$. Hence, $(l_r, \tau_q)$ is a quasi-locally convex topological vector space.

More precisely, from (5.5), the topology $\tau_q$ on $l_r$ defined by the mapping $q$ is the coarsest topology on $X$, for which, the following mappings

$$q_\lambda(\cdot - \{u_i\}): l_r \to \mathbb{R}^+, \text{ for any } \{u_i\} \in l_r \text{ and for every } \lambda = 1, 2, \ldots, \text{ are } \tau_q\text{-continuous.} \tag{5.6}$$

We next prove (5.6). From the triangle inequality in $l_r$, for $0 < r < 1$, we have

$$q(\{t_i\} + \{s_i\}) \leq q(\{t_i\}) + q(\{s_i\}), \text{ for } \{t_i\}, \{s_i\} \in l_r,$$

and

$$q_\lambda(\{t_i\} + \{s_i\}) \leq q_\lambda(\{t_i\}) + q_\lambda(\{s_i\}), \text{ for } \lambda = 1, 2, \ldots.$$

For any $\{u_i\} \in l_r$, it follows

$$q_\lambda(\{t_i\} - \{u_i\}) = q_\lambda(\{t_i\} - \{s_i\} + \{s_i\} - \{u_i\}) \leq q_\lambda(\{t_i\} - \{s_i\}) + q_\lambda(\{s_i\} - \{u_i\}).$$

It implies

$$|q_\lambda(\{t_i\} - \{u_i\}) - q_\lambda(\{s_i\} - \{u_i\})| \leq q_\lambda(\{t_i\} - \{s_i\}) \leq q(\{t_i\} - \{s_i\}). \tag{5.7}$$

Then (5.6) follows from (5.7) immediately. Therefore, we proved that ($l_r$, $\tau_q$) is quasi-locally convex.

**References**

[1] L.E. J. Brouwer, Uber Abbildung von Mannigfaltigkeiten, Math. Ann. 7l (1912) 97–1l5.

[2] F. F. Bonsall, *Lectures on some fixed point theorems of functional analysis*, Bombay (1962).

[3] Butsan, T., Dhompongsa, S., Fupinwong, W., Schuader's conjecture on convex metric spaces, J. Nonlinear Convex Anal. 11, no. 3, (2010) 527–535.

[4] R. Cauty in Solution du problème de point fixe de Schauder, Fund. Math. 170 (2001) 231–246]

[5] K. Fan, A generalization of Tychonoff's fixed point theorem, *Math. Ann*., **142** (1961) 305–310.

[6] Olga Hadžić, A fixed point theorem in topological vector spaces, Review of Research Faculty of Science, University of Novi Sad, Volume 10 (1980).

[7] Won Kyu Kim**,** A fixed point theorem in a Hausdorff topological vector space, Comment. Math. Univ. Carolin. 36, 1 (1995)33–38.

[8] Sehie Park**,** The KKM principle implies many fixed point theorems, Topology and its Applications 135 (2004) 197–206

[9] S. Park, Recent applications of Fan-KKM theorem, Lecture notes in Math. Anal. Kyoto University, **1841**, 58–68 (2013).

[10] J. Schauder, *Der Fixpunktsatz in Funktionalräumen*, Studia Math. 2 (1930) 171–180

[11] Joel H. Shapiro, The Schauder Fixed-Point Theorem, An Infinite Dimensional Brouwer Theorem, A Fixed-Point Farrago pp 75–81

[12] A. Tychonoff, *Ein Fixpunktsatz*, Mathematische Annalen 111 (1935) 767–776

## Appendix

A proof that $(l_r, \tau_q)$ is not locally convex.

Assume, on the contrary, that $(l_r, \tau_q)$ is locally convex. By the second version of the definition of locally convex topological vector spaces, we suppose that $(l_r, \tau_q)$ has a base at the origin consisting of convex open sets, denoted by $\{C_\gamma\}$. On the other hand, the $q$-topology on $(l_r, \tau_q)$ insures that the origin of $l_r$ has a base consisting of the following family of open balls

$$B(\delta^r) = \{\{t_i\} \in l_r : \sum_{i=1}^{\infty} |t_i|^r < \delta^r\}, \delta > 0.$$

So, for any fixed $\delta > 0$, there are $C_0 \in \{C_\gamma\}$ and $\alpha \in (0, \delta)$ such that

$$B(\alpha^r) \subseteq C_0 \subseteq B(\delta^r).$$

Then, for any $\beta \in (0, \alpha)$, we have that, for any $\{t_i\} \in l_r$,

$$q(\{t_i\}) = \beta^r \text{ implies } \{t_i\} \in B(\alpha^r) \subseteq C_0 \subseteq B(\delta^r). \tag{5.8}$$

For any positive number $j \geq 1$, with respect to $\beta$, we define

$$\Delta_j = (0, 0, \ldots, \beta, 0, 0, \ldots), \text{ for } j = 1, 2, \ldots,$$

where $\beta$ is the $j^{\text{th}}$-coordinate. It follows that $q(\Delta_j) = \beta^r$; and therefore $\Delta_j \in B(\beta^r) \subseteq C_0$, for $j = 1, 2, \ldots$. Since $C_0$ is convex, then, for any positive integer $n$, from (5.8), we have

$$\sum_{j=1}^{n} \frac{1}{n} \Delta_j \in C_0 \subseteq B(\delta^r), \tag{5.9}$$

and

$$q\left(\sum_{j=1}^{n} \frac{1}{n} \Delta_j\right) = \sum_{j=1}^{n} \left(\frac{1}{n}\beta\right)^r = n^{1-r}\beta^r. \tag{5.10}$$

From (5.9) and (5.10), we have

$$\delta^r \geq q\left(\sum_{j=1}^{n} \frac{1}{n} \Delta_j\right) = n^{1-r}\beta^r. \tag{5.11}$$

Since $\beta > 0$, $0 < r < 1$, and $n$ can be arbitrary large, (5.11) creates a contradiction. Hence $(l_r, \tau_q)$ is not locally convex.